\documentclass[fleqn, a4paper, oneside, doublespacing]{imsart}

\usepackage{amsthm,amsmath,amssymb}
\usepackage[authoryear]{natbib}
\bibpunct{(}{)}{;}{a}{,}{,}
\RequirePackage[OT1]{fontenc}

\startlocaldefs

\newtheorem{theorem}{Theorem}
\newtheorem{corollary}[theorem]{Corollary}
\newtheorem{lemma}[theorem]{Lemma}
\newtheorem{proposition}[theorem]{Proposition}

\theoremstyle{remark}
\newtheorem{remark}[theorem]{Remark}

\theoremstyle{definition}
\newtheorem{definition}{Definition}

\def\R{\mathbb{R}}

\endlocaldefs

\begin{document}

\begin{frontmatter}

\title{Data-driven efficient score tests for deconvolution problems.}
\runtitle{Score tests for deconvolution.}

\thankstext{t2}{Financial support of the Deutsche
Forschungsgemeinschaft GK 1023 "Identifikation in Mathematischen
Modellen" is gratefully acknowledged.}

\begin{aug}
\author{\fnms{Langovoy}
\snm{Mikhail}\thanksref{t2}\ead[label=e1]{langovoy@math.uni-goettingen.de}}

\affiliation{Institute for Mathematical Stochastics,\\
         Georg-August-University at G\"{o}ttingen,\\
         Germany.}

\address{Mikhail Langovoy, Institute for Mathematical Stochastics, \\
Georg-August-University at G\"{o}ttingen, Maschmuellenweg 8/10,
\\
37077 Goettingen, Germany\\
\printead{e1}\\}

\runauthor{M. Langovoy}
\end{aug}

\begin{abstract}
We consider testing statistical hypotheses about densities of
signals in deconvolution models. A new approach to this problem is
proposed. We constructed score tests for the deconvolution with
the known noise density and efficient score tests for the case of
unknown density. The tests are incorporated with model selection
rules to choose reasonable model dimensions automatically by the
data. Consistency of the tests is proved.
\end{abstract}



\end{frontmatter}

\section{Introduction}\label{SimpleDeconvolutionSection1}

Constructing good tests for statistical hypotheses is an essential
problem of statistics. There are two main approaches to
constructing test statistics. In the first approach, roughly
speaking, some measure of distance between the theoretical and the
corresponding empirical distributions is proposed as the test
statistic. Classical examples of this approach are the Cramer-von
Mises and the Kolmogorov-Smirnov statistics. Although, these tests
works and are capable of giving very good results, but each of
these tests is asymptotically optimal only in a finite number of
directions of alternatives to a null hypothesis (see
\cite{MR1335235}).

Nowadays, there is an increasing interest to the second approach
of constructing test statistics. The idea of this approach is to
construct tests in such a way that the tests would be
asymptotically optimal. Test statistics constructed following this
approach are often called (efficient) score test statistics. The
pioneer of this approach was \cite{Neyman1937} and then many other
works followed: \cite{MR0112201}, \cite{MR0370837},
\cite{MR1226715}, \cite{MR0084918}, \cite{MR1294744}. This
approach is also closely related to the theory of efficient
(adaptive) estimation - \cite{MR1245941}, \cite{MR620321}. Score
tests are asymptotically optimal in the sense of intermediate
efficiency in an infinite number of directions of alternatives
(see \cite{MR1421157}) and show good overall performance in
practice (\cite{MR1370299}, \cite{MR1488856}).

We described the situation in classical hypothesis testing, i.e.,
testing hypotheses about random variables $X_1, \ldots , X_n,$
whose values are directly observable. But, it is important from
practical point of view to be able to construct tests for
situations where $X_1, \ldots , X_n$ are corrupted or can only be
observed with an additional noise term. These kind of problems are
termed \emph{statistical inverse problems}. The most well-known
example here is the deconvolution problem. This problem appears
when one has noisy signals or measurements: in physics,
seismology, optics and imaging, engineering. It is a building
block for many complicated statistical inverse problems.

Due to importance of the deconvolution problem, testing
statistical hypotheses related to this problem has been widely
studied in the literature. But, to our knowledge, all the proposed
tests were based on some kind of distance (usually a $L_2-$type
distance) between the theoretical density function and the
empirical estimate of the density (see, for example,
\cite{MR0348906}, \cite{MR1979392}, \cite{HBM2007}). Thus, only
the first approach described above was implemented for the
deconvolution problem.

In this paper, we treat the deconvolution problem with the second
approach. We construct efficient score tests for the problem. From
classical hypothesis testing, it was shown that for applications
of efficient score tests, it is important to select the right
number of components in the test statistic (see \cite{MR1226715},
\cite{MR1248222}, \cite{MR1370299}, \cite{MR1395735}). Thus, we
provide corresponding refinement of our tests. Following the
solution proposed in \cite{MR1964450}, we make our tests
data-driven, i.e., the tests are capable to choose a reasonable
number of components in the test statistics automatically by the
data.

In Section \ref{SimpleDeconvolutionSection2}, we formulate the
simple deconvolution problem. In Section
\ref{SimpleDeconvolutionSection4}, we construct the score tests
for the parametric deconvolution hypothesis. In Section
\ref{SimpleDeconvolutionSection5}, we prove consistency of our
tests against nonparametric alternatives. In Section
\ref{CompositeDeconvolutionSection2}, we turn to the deconvolution
with an unknown error density. We derive the efficient scores for
the composite parametric deconvolution hypothesis in Section
\ref{CompositeDeconvolutionSection3}. In Section
\ref{CompositeDeconvolutionSection4}, we construct the efficient
score tests for this case. In Section
\ref{CompositeDeconvolutionSection5}, we make our tests
data-driven. In Section \ref{CompositeDeconvolutionSection6}, we
prove consistency of the tests against nonparametric alternatives.
Additionally, in Sections \ref{SimpleDeconvolutionSection5} and
\ref{CompositeDeconvolutionSection6}, we explicitly characterize
the class of nonparametric alternatives such that our tests are
inconsistent and therefore shouldn't be used for testing against
the alternatives from this class. Some simple examples of
applications of the theory are also presented in this paper.

\section{Notation and basic assumptions}\label{SimpleDeconvolutionSection2}

The problem of testing whether i.i.d. real-valued random variables
$X_1, \ldots , X_n$ are distributed according to a given density
$f$ is classical in statistics. We consider a more difficult
problem, namely the case when $X_i$ can only be observed with an
additional noise term, i.e., instead of $X_i$ one observes $Y_i,$
where
$$Y_i = X_i + \varepsilon_i ,$$ and $\varepsilon_i'$s are i.i.d. with a known density $h$ with respect
to the Lebesgue measure $\lambda;$ also $X_i$ and $\varepsilon_i$
are independent for each $i$ and $E\,\varepsilon_i =0,\,
 0 < E\,{\varepsilon}^2 < \infty .$ For brevity of notation say that $X_i,$ $Y_i,$
$\varepsilon_i$ have the same distribution as random variables
$X,$ $Y,$ $\varepsilon$ correspondingly. Assume that $X$ has a
density with respect to $\lambda.$

Our null hypothesis $H_0$ is the simple hypothesis that $X$ has a
known density $f_0$ with respect to $\lambda.$ The most general
possible nonparametric alternative hypothesis $H_A$ is that $f
\neq f_0.$ Since this class of alternatives is too broad, first we
would be concerned with a special class of submodels of the model
described above. In this paper we will at first assume that all
possible alternatives from $H_A$ belong to some parametric family.
Then we will propose a test that is expected to be asymptotically
optimal (in some sense) against the alternatives from this
parametric family. However, we will prove that our test is
consistent also against other alternatives even if they do not
belong to the initial parametric family. The test is therefore
applicable in many nonparametric problems. Moreover, the test is
expected to be asymptotically optimal (in some sense) for testing
against an infinite number of directions of nonparametric
alternatives (see \cite{MR1421157}). This is the general plan for
our construction. \vspace{6pt}

\section{Score test for simple deconvolution}\label{SimpleDeconvolutionSection4}

Suppose that all possible densities of $X$ belong to some
parametric family $\{f_{\theta}\},$ where $\theta$ is a
$k-$dimensional Euclidean parameter, $\Theta \in {\R}^{k}$ is a
parameter set. Then all the possible densities $q\,(y\,;\theta)$
of $Y$ have in such model the form
\begin{equation}\label{SimpleDeconvolution13}
q\,(y\,;\theta)= \int_{\R}\,f_{\theta}(s)\,h(\,y-s)\,ds\,.
\end{equation}

\noindent The \emph{score function} $\dot{l}$ is defined as

\begin{equation}\label{SimpleDeconvolution14}
\dot{l}\,(y\,;\theta)= \,
\frac{{\bigr(q\,(\theta)\bigr)}_{\theta}^{'} } {\,q\,(\theta)}
\,1_{[q\,(\theta)>0]}\,,
\end{equation}

\noindent where $q\,(\theta):=q\,(y;\theta)$ and
$l\,(\theta):=l\,(y;\theta)$ for brevity. The \emph{Fisher
information matrix} of parameter $\theta$ is defined as

\begin{equation}\label{SimpleDeconvolution15}
I(\theta)= \,\int_{\R} \dot{l}\,(y\,;\theta)
\dot{l}^{T}\,(y\,;\theta)\,dQ_{\theta}(y)\,.
\end{equation}

\begin{definition}\label{SimpleDeconvolutionDef2}
Call our problem a \emph{regular deconvolution problem} if

$\quad$

$\langle \mathbf{B1} \rangle \quad\quad\mbox{for all}\,\, \theta
\in \Theta \quad q\,(y;\theta) \,\, \mbox{is continuously
differentiable in}\,\, \theta$

$\quad\quad\quad\quad\,\mbox{for}\,\, \lambda-\mbox{almost
all}\,\, y \,\,\mbox{with gradient}\,\, \dot{q}\,(\theta)$

$\quad$

$\langle \mathbf{B2} \rangle \quad\quad \bigr|\dot{l}\,(\theta)
\bigr| \in L_2 (\R, Q_\theta)\,\,\,\mbox{for all}\,\, \theta \in
\Theta$

$\quad$

$\langle \mathbf{B3} \rangle \quad\quad I(\theta) \quad\mbox{is
nonsingular for all}\,\, \theta \in \Theta \,\, \mbox{and
continuous in}\,\, \theta\,.$
\end{definition}

\noindent If $\theta$ is a true parameter value, call such model
$GM_k(\theta)$ and denote by $Q_\theta$ the probability
distribution function and by $E_\theta$ the expectation
corresponding to the density $q\,( \cdot ;\theta).$

If conditions $\langle B1 \rangle - \langle B3 \rangle$ holds,
then by Proposition 1, p.13 of \cite{MR1245941} we calculate for
all $y \in supp\,\,q\,( \cdot ;\theta)$

\begin{equation}\label{SimpleDeconvolution16}
\dot{l}\,(\theta)=\dot{l}\,(y\,;\theta)= \,
\frac{{\bigr(q\,(y\,;\theta)\bigr)}_{\theta}^{'} }
{\,q\,(y\,;\theta)} \,= \frac{ \frac{ \,\partial
}{\,\partial\theta}\,\int_{\R}
\,f_{\theta}(s)\,h(\,y-s)\,ds}{\int_{\R}\,f_{\theta}(s)\,h(\,y-s)\,ds}\,.
\end{equation}

\noindent Then for $y \in supp\,\,q\,( \cdot ;\theta)$ the
\emph{efficient score vector} for testing $H_0:\,\theta =0$ is

\begin{equation}\label{SimpleDeconvolution17}
l^\ast (y) := \dot{l}\,(y\,;0)\,= \frac{ \frac{ \,\partial
}{\,\partial\theta}\,\Bigr(\int_{\R}\,f_{\theta}(s)\,h(\,y-s)\,ds\Bigr)\Bigr|_{\theta=0}}
{\int_{\R}\,f_0(s)\,h(\,y-s)\,ds}\,.
\end{equation}

Set
\begin{equation}\label{SimpleDeconvolution18}
L= {\{ E_0 {[l^\ast (Y)]}^{T} l^\ast (Y) \}}^{-1}
\end{equation}
\noindent and
\begin{equation}\label{SimpleDeconvolution19}
U_k =\biggr\{ \frac{1}{\sqrt{n}}\, \sum_{j=1}^{n} l^\ast (Y_j)
\biggr\} \;L\; {\biggr\{ \frac{1}{\sqrt{n}}\, \sum_{j=1}^{n}
l^\ast (Y_j) \biggr\} }^{T} \,.
\end{equation}

\begin{theorem}\label{SimpleDeconvolutionTh3}
For the regular deconvolution problem the efficient score vector
$l^\ast$ for testing $\theta=0$ in $GM_k(\theta)$ is given for all
$x \in \R$ by (\ref{SimpleDeconvolution17}). Moreover, under $H_0
: \theta = 0$ we have $U_k \rightarrow_d \;{\chi}^{2}_{k}$ as $n
\rightarrow \infty.$
\end{theorem}

We construct the test based on the test statistic $U_k$ as
follows: the null hypothesis $H_0$ is rejected if the value of
$U_k$ exceeds standard critical points for
${\chi}^{2}_{k}-$distribution. Note that we do not need to
estimate the scores $l^{\ast}.$

\begin{corollary}\label{SimpleDeconvolutionCor4}
If the deconvolution problem is regular and $f_{\theta}(\cdot)$ is
differentiable in $\theta$ for all $\theta \in \Theta,$ then the
conclusions of Theorem \ref{SimpleDeconvolutionTh3} are valid and
the efficient score vector for testing $H_0 : \theta = 0$ can be
calculated by the formula

\begin{equation}\label{SimpleDeconvolution20}
l^\ast (y)\,=
\frac{\int_{\R}\bigr(\frac{\,\partial}{\,\partial\theta}\,
f_{\theta}(s)\,\bigr)\bigr|_{\theta=0}\, h(\,y-s)\,ds}
{\int_{\R}\,f_{\theta}(s)\,h(\,y-s)\,ds}\,.
\end{equation}
\end{corollary}


\noindent {\bf Example 1.} Consider one important special case.
Assume that each submodel of interest is given by the following
restriction: all possible densities $f$ of $X$ belong to a
parametric exponential family, i.e., $f=f_\theta$ for some
$\theta,$ where

\begin{equation}\label{SimpleDeconvolution1}
f_{\theta}(x)\,=\,f_0(x)\,b(\theta)\exp (\theta \circ u(x))\,,
\end{equation}

\noindent where the symbol $\circ$ denotes the inner product in
${\R}^k,$ $u(x) = (u_1(x), \ldots, u_k(x))$ is a vector of known
Lebesgue measurable functions, $b(\theta)$ is the normalizing
factor and $\theta\in \Theta \subseteq {\R}^k.$ We assume that the
standard regularity assumptions on exponential families (see
\cite{MR489333}) are satisfied. All the possible densities
$q\,(y\,;\theta)$ of $Y$ have in such model the form

\begin{equation}\label{SimpleDeconvolution4}
q\,(y\,;\theta)= \int_{\R}\,f_0(s)\,b(\theta)\exp (\theta \circ
u(s))\,\,h(\,y-s)\,ds\,.
\end{equation}

\noindent These densities no longer need to form an exponential
family. If we assume, for example, that $h>0$ $\quad \lambda
-\mbox{almost everywhere on} \,\,\,\mathbb{R}$ and the functions
$f_0,$ $h,$ $u_1, \ldots, u_k$ are bounded and $\lambda
-\mbox{measurable}$ and that there exists an open subset $\Theta_1
\subseteq \Theta$ such that $\bigr| \dot{l}\,(y\,;\theta) \bigr|
\in L_2 (Q_\theta)$ and the Fisher information matrix $I(\Theta)$
is nonsingular and continuous in $\theta ,$ then conditions
$\langle B1 \rangle - \langle B3 \rangle$ are satisfied for this
problem and the previous results are applicable. The score vector
for the problem is
\begin{equation}\label{SimpleDeconvolution7}
l^\ast (y) = \,\frac{\int_{\R}\, u(s)\,
f_0(s)\,h(\,y-s)\,ds\,}{\int_{\R}\,f_0(s)\,h(\,y-s)\,ds}\,-\,
\int_{\R}\,u(s)\,f_0(s)\,ds\,.
\end{equation}













\noindent In other words, if we denote by $\ast$ the standard
convolution of functions,
\begin{equation}\label{SimpleDeconvolution8}
l^\ast (y) = \, \frac{(u f_0)\ast h}{f_0 \ast h}(y)\,-\, E_0
u(X)\,.
\end{equation}
\noindent Let $L$ be defined by (\ref{SimpleDeconvolution18}) and


\begin{equation}\label{SimpleDeconvolution10}
V_k =\biggr\{ \frac{1}{\sqrt{n}}\, \sum_{j=1}^{n} l^\ast (Y_j)
\biggr\} \;L\; {\biggr\{ \frac{1}{\sqrt{n}}\, \sum_{j=1}^{n}
l^\ast (Y_j) \biggr\} }^{T} \,.
\end{equation}

\noindent This is the score test statistic designed to be
asymptotically optimal for testing $H_0$ against the alternatives
from the exponential family (\ref{SimpleDeconvolution1}). Its
asymptotic distribution under the null hypothesis $H_0$ is given
by Theorem \ref{SimpleDeconvolutionTh3}.

\vspace{6pt}

\section{Selection rule}\label{SimpleDeconvolutionSubsection3.1}

For the use of score tests in classical hypotheses testing it was
shown (see the Introduction) that it is important to select the
right dimension $k$ of the space of possible alternatives.
Incorrect choice of the model dimension can substantially decrease
the power of a test. In Section \ref{SimpleDeconvolutionSection5}
we give a theoretical explanation of this fact for the case of
deconvolution. The possible solution of this problem is to
incorporate the test statistic of interest by some procedure
(called a selection rule) that chooses a reasonable dimension of
the model automatically by the data. See \cite{MR1964450} for an
extensive discussion and practical examples. In this section we
implement this idea for testing the deconvolution hypothesis.
First we give a definition of selection rule, generalizing ideas
from \cite{MR2259976}.

Denote by $M_k (\theta)$ the model described in Section
\ref{SimpleDeconvolutionSection4} such that the true parameter
$\theta$ belongs to the parameter set, say $\Theta_k ,$ and $\dim
\Theta_k = k.$ By a \emph{nested family} of submodels $M_k
(\theta)$ for $k=1,2, \ldots$ we mean a sequence of these models
such that for their parameter sets it holds that $\Theta_1$
$\subset$ $\Theta_2$ $\subset \ldots .$

\begin{definition}\label{SimpleDeconvolutionDef1}
Consider a nested family of submodels $M_k (\theta)$ for $k=1,
\ldots , d,$ where $d$ is fixed but otherwise arbitrary. Choose a
function $\pi (\cdot,\cdot): \mathbb{N} \times \mathbb{N}
\rightarrow \R,$ where $\mathbb{N}$ is the set of natural numbers.
Assume that $\pi(1,n) <\pi(2,n)< \ldots <\pi(d,n)$ for all $n$ and
$\pi(j,n)-\pi(1,n) \rightarrow \infty$ as $n \rightarrow \infty$
for every $j= 2, \ldots ,d.$ Call $\pi(j,n)$ a \emph{penalty
attributed to the j-th model $M_j(\theta)$ and the sample size n.}
Then a \emph{selection rule} $S$ for the test statistic $U_k$ is
an integer-valued random variable satisfying the condition

\begin{equation}\label{SimpleDeconvolution11}
S=\min\bigr\{ k:\,1\leq k \leq d;\, U_k - \pi(k,n) \geq U_j
-\pi(j,n),\, j= 1, \ldots ,d\bigr\}\,.
\end{equation}

\noindent We call $U_S$ a \emph{data-driven efficient score test
statistic} for testing validity of the initial model.
\end{definition}

From Theorem \ref{SimpleDeconvolutionTh2} below it follows that
for our problem (as well as in the classical case, see
\cite{MR1964450}) many possible penalties lead to consistent
tests. So the choice of the penalty should be dictated by external
practical considerations. Our simulation study is not so vast to
recommend the most practically suitable penalty for the
deconvolution problem. Possible choices are, for example,
Schwarz's penalty $\pi(j,n)=j \, \log n,$ or Akaike's penalty
$\pi(j,n)=j.$

Denote by $P^n_0$ the probability measure corresponding to the
case when $X_1, \ldots , X_n$ all have the density $f_0.$ For
simplicity of notation we will further sometimes omit index "n"
and write simply $P_0.$ The main result about the asymptotic null
distribution of $U_S$ is the following

\begin{theorem}\label{SimpleDeconvolutionTh2}
Suppose that assumptions $\langle B1 \rangle-\langle B3 \rangle$
holds. Then under the null hypothesis $H_0$ it holds that $P^n_0
(S>1) \rightarrow 0$ and $U_S \rightarrow_d \,{\chi}^{2}_{1}$ as
$n \rightarrow \infty .$
\end{theorem}

\begin{remark}\label{SimpleDeconvolutionRemark1} The selection rule
$S$ can be modified in order to make it possible to choose not
only models of dimension less than some fixed $d$ but to allow
arbitrary large dimensions of $M_k(\theta)$ as $n$ grows to
infinity. In this case an analogue of Theorem
\ref{SimpleDeconvolutionTh2} still holds, but the proof becomes
more technical and one should take care about the possible rates
of growth of the model dimension. Though, one can argue that even
$d=10$ is often enough for practical purposes (see
\cite{MR1370299}).

\end{remark}
\vspace{6pt}

\section{Consistency of tests}\label{SimpleDeconvolutionSection5}

Let $F$ be a true distribution function of $X.$ Here $F$ is
\emph{not} necessarily parametric and possibly doesn't have a
density with respect to $\lambda .$ Let us choose for every $k
\leq d$ an auxiliary parametric family $\{f_{\theta}\},\,$ $\theta
\in \Theta \subseteq \mathbb{R}^k$ such that $f_0$ from this
family coincides with $f_0$ from the null hypothesis $H_0.$
Suppose that the chosen family $\{f_{\theta}\}$ gives us the
regular deconvolution problem in the sense of Definition
\ref{SimpleDeconvolutionDef2}. Then one is able to construct the
score test statistic $U_k$ defined by
(\ref{SimpleDeconvolution19}) despite the fact that the true $F$
possibly has no relation to the chosen $\{f_{\theta}\}.$ One can
use the exponential family from Example 1 as $\{f_{\theta}\},$ or
some other parametric family whatever is convenient. This is our
goal in this section to determine under what conditions thus build
$U_k$ will be consistent for testing against $F.$

Suppose that the following condition holds

$\quad$

$\langle \mathbf{D1} \rangle \quad\quad\mbox{there exists an
integer}\,\, K \geq 1 \,\,\mbox{such that}\,\, K \leq d
\,\,\mbox{and}$

$\quad\quad\quad\quad E_F\,l_1^{\ast}=0, \ldots ,
E_F\,l_{K-1}^{\ast}=0,\,\,E_F\,l_K^{\ast}= C_K \neq 0\,,$

$\quad$

\noindent where $l_i^{\ast}$ is the $i-$th coordinate function of
$l^{\ast}$ and $l^{\ast}$ is defined by
(\ref{SimpleDeconvolution17}), $d$ is the maximal possible
dimension of our model as in Definition
\ref{SimpleDeconvolutionDef1} of Section
\ref{SimpleDeconvolutionSubsection3.1}, and $E_F$ denotes the
mathematical expectation with respect to $F \ast h.$

Condition $\langle D1 \rangle$ is a weak analog of nondegeneracy:
if for all $k$ $\langle D1 \rangle$ fails, then $F$ is orthogonal
to the whole system ${\{l_i^{\ast} \}}_{i=1}^{\infty},$ and if
this system is complete, then $F$ is degenerate. Also $\langle D1
\rangle$ is related to the identifiability of the model (see the
beginning of Section \ref{CompositeDeconvolutionSection6} for more
details).

We start with investigation of consistency of $U_k,$ where $k$ is
some fixed number, $1 \leq k \leq d.$ The following result shows
why it is important to choose the right dimension of the model.

\begin{proposition}\label{SimpleDeconvolutionProp5}
Let $\langle D1 \rangle$ holds. Then for all $1 \leq k \leq K-1,$
if $F$ is the true distribution function of $X,$ then $U_k
\rightarrow_d \,{\chi}^{2}_{k}$ as $n \rightarrow \infty\,.$
\end{proposition}

\noindent This result and Theorem \ref{SimpleDeconvolutionTh3}
show that if the dimension of the model is too small, then the
test doesn't work since it doesn't distinguish between $F$ and
$f_0.$

\begin{proposition}\label{SimpleDeconvolutionProp6}
Let $\langle D1 \rangle$ holds. Then for $k \geq K,$ if $F$ is the
true distribution function of $X,$ then $U_k \rightarrow \,\infty$
in probability as $n \rightarrow \infty\,.$
\end{proposition}

\noindent Now we turn to the data-driven statistic $U_S.$ Suppose
that the selection rule $S$ is defined as in Section
\ref{SimpleDeconvolutionSubsection3.1}. Assume that

$\quad$

$\langle \mathbf{S1} \rangle \quad\quad\mbox{for every fixed }\,\,
k \geq 1 \,\,\mbox{it holds that}\,\, \pi (k,n)= o(n)
\,\,\mbox{as}\,\, n \rightarrow \infty\,.$

$\quad$

\noindent Denote by $P_F$ the probability measure corresponding to
the case when $X_1,$ $\ldots ,$ $X_n$ all have the distribution
$F.$ Consider consistency of the "adaptive" test based on $U_S.$

\begin{proposition}\label{SimpleDeconvolutionProp7}
Let $\langle D1 \rangle$ and $\langle S1 \rangle$ holds. If $F$ is
the true distribution function of $X,$ then $P_F \,(S \geq K)
\rightarrow 1$ and $U_S \rightarrow \,\infty$ as $n \rightarrow
\infty\,.$
\end{proposition}

\noindent The main result of this section is the following

\begin{theorem}\label{SimpleDeconvolutionTh8}
$\quad$

\begin{enumerate}
\item \noindent The test based on $U_k$ is consistent for testing
against all alternative distributions $F$ such that $\langle D1
\rangle$ is satisfied with $K \leq k$

\item \noindent The test based on $U_k$ is inconsistent for
testing against all alternative distributions $F$ such that
$\langle D1 \rangle$ is satisfied with $K > k$

\item \noindent If the selection rule $S$ satisfies $\langle S1
\rangle,$ then test based on $U_S$ is consistent against all
alternative distributions $F$ such that $\langle D1 \rangle$ is
satisfied with some $K.$
\end{enumerate}
\end{theorem}

\section{Composite deconvolution}\label{CompositeDeconvolutionSection2}

In the previous sections, we treated the simplest case of the
deconvolution problem. The next sections are devoted to the more
realistic case of unknown error density. Our main ideas and
constructions will be similar to the ones for the simple case. Our
goal is to modify the technics and constructions from the simple
hypothesis case in order to apply them in the new situation. In
order to do this we will have to impose on our new model
additional regularity assumptions of uniformity. These assumptions
are quite standard in statistics. They are a necessary payment for
our ability to keep simple and general constructions for the more
complicated problem. We will have to modify the scores we used in
the simple case. The modification we will use is called efficient
scores.

Despite of all the changes, we will still be able to build a
selection rule for the new problem. We will need a new and
modified definition of the selection rule. Big part of the new
model uniformity assumptions will be needed not to build an
efficient score test, but to make such test data-driven  (see
section \ref{CompositeDeconvolutionSection5}).


Consider the situation described in the first paragraph of Section
\ref{SimpleDeconvolutionSection2}, but with the following
complication introduced. Suppose further on that the density $h$
of $\varepsilon$ is \emph{unknown}.

Then the most general possible null hypothesis $H_0$ in this setup
is that $f=f_0$ and the error $\varepsilon$ has expectation 0 and
finite variance. The most general alternative hypothesis $H_A$ is
that $f \neq f_0.$ Since both $H_0$ and $H_A$ are in this case too
broad, we would first consider a special class of submodels of the
model described above. At first we assume that all possible
densities $f$ of $X$ belong to some specific and preassigned
parametric family $\{ f_\theta \},$ i.e., $f=f_\theta$ for some
$\theta$ and $\theta$ is a $k-$dimensional Euclidian parameter and
$\Theta \subseteq R^k$ is a parameter set for $\theta.$ Our
starting assumption about the density of the error $\varepsilon$
will be that $h$ belongs to some specific parametric family $\{
h_\eta \},$ where $\eta \in \Lambda$ and $\Lambda \subseteq R^m$
is a parameter set. Thus, $\eta$ is a nuisance parameter. The null
hypothesis $H_0$ is the following composite hypothesis: $X$ has
particular density $f_0$ with respect to $\lambda.$

Then we will propose a test that is expected to be asymptotically
optimal (in some sense) for testing in this parametric situation.
After that we will prove that our test is consistent also against
a wide class of nonparametric alternatives. Moreover, the test is
expected to be asymptotically optimal (in some sense) for testing
against an infinite number of directions of nonparametric
alternatives. This is essentially the same plan as for the simple
case.

If $(\theta, \eta)$ is a true parameter value, we call such
submodel $M_{k,m}(\theta, \eta).$ Denote in this case the density
of $Y$ by $g(\cdot;(\theta, \eta))$ and the corresponding
expectation by $E_{(\theta, \eta)}.$ Let the null hypothesis $H_0$
be $\theta=\theta_0,$ where it is assumed that $\theta_0 \in
\Theta.$ Then the alternative hypothesis $\theta \neq \theta_0$ is
a parametric subset of the original general and nonparametric
alternative hypothesis $H_A.$\vspace{6pt}

\section{Efficient scores}\label{CompositeDeconvolutionSection3}

All possible densities $g\,(y\,;(\theta,\eta))$ of $Y$ have in our
model the form
\begin{equation}\label{CompositeDeconvolution1}
g\,(y\,;(\theta,\eta))=
\int_{\R}\,f_{\theta}(s)\,h_\eta(\,y-s)\,ds\,.
\end{equation}

\noindent It is not always possible to identify $\theta$ or/and
$\eta$ in this model. Since we are concerned with testing
hypotheses and not with estimation of parameters, it is not
necessary for us to impose a restrictive assumption of
identifiability on the model. We will need only a (weaker)
consistency condition to build a sensible test (see Section
\ref{CompositeDeconvolutionSection6}).

The \emph{score function} for $(\theta, \eta)$ at $(\theta_0,
\eta_0)$ is defined as (see \cite{MR1245941}, p.28):

\begin{equation}\label{CompositeDeconvolution2}
\dot{l}_{\theta_0, \eta_0}(y)\,=\,\bigr(\,\dot{l}_{\theta_0}(y),\,
\dot{l}_{\eta_0}(y)\,\bigr)\,,
\end{equation}

\noindent where $\dot{l}_{\theta_0}$ is the score function for
$\theta$ at $\theta_0$ and $\dot{l}_{\eta_0}$ is the score
function for $\eta$ at $\eta_0,$ i.e.

\begin{equation}\label{CompositeDeconvolution3}
\dot{l}_{\theta_0}(y)= \, \frac{ \frac{\,\partial
}{\,\partial\theta}\,\bigr( \,g\,(y\,;(\theta,\eta_0))
\bigr)\bigr|_{\theta=\theta_0} } {\,g\,(y\,;(\theta_0,\eta_0))}
\,1_{[y:\,g\,(y\,;(\theta_0,\eta_0))>0]}\,
\end{equation}

$$=\,\frac{ \frac{\,\partial }{\,\partial\theta}\,\Bigr(
\int_{\R}\,f_{\theta}(s)\,h_{\eta_0}(\,y-s)\,ds
\Bigr)\Bigr|_{\theta=\theta_0} }
{\int_{\R}\,f_{\theta_0}(s)\,h_{\eta_0}(\,y-s)\,ds }
\,1_{[y:\,g\,(y\,;(\theta_0,\eta_0))>0]}\,,$$

\begin{equation}\label{CompositeDeconvolution4}
\dot{l}_{\eta_0}(y)= \, \frac{ \frac{\,\partial
}{\,\partial\eta}\,\bigr( \,g\,(y\,;(\theta_0,\eta))
\bigr)\bigr|_{\eta=\eta_0} } {\,g\,(y\,;(\theta_0,\eta_0))}
\,1_{[y:\,g\,(y\,;(\theta_0,\eta_0))>0]}\,
\end{equation}

$$=\,\frac{ \frac{\,\partial }{\,\partial\eta}\,\Bigr(
\int_{\R}\,f_{\theta_0}(s)\,h_{\eta}(\,y-s)\,ds
\Bigr)\Bigr|_{\eta=\eta_0} }
{\int_{\R}\,f_{\theta_0}(s)\,h_{\eta_0}(\,y-s)\,ds }
\,1_{[y:\,g\,(y\,;(\theta_0,\eta_0))>0]}\,.$$

\noindent The \emph{Fisher information matrix} of parameter
$(\theta, \eta)$ is defined as

\begin{equation}\label{CompositeDeconvolution5}
I(\theta, \eta)\,= \,\int_{\R} {\dot{l}_{\theta,
\eta}}^T(y)\,\dot{l}_{\theta, \eta}(y)\, \,dG_{\theta, \eta}(y)\,,
\end{equation}

\noindent where $G_{\theta, \eta}(y)$ is the probability measure
corresponding to the density $g\,(y\,;(\theta,\eta)).$ The symbol
'T' denotes the transposition and all vectors are supposed to be
row ones.

We assume that $M_{k,m}(\theta, \eta)$ is a regular parametric
model in the sense of the following definition.

\begin{definition}\label{CompositeDeconvolutionDefinition1}
Call our problem a \emph{regular deconvolution problem} if

$\quad$

$\langle \mathbf{A1} \rangle \quad\quad\mbox{for all}\,\, (\theta,
\eta) \in \Theta \times \Lambda \quad g\,(y\,;(\theta,\eta)) \,\,
\mbox{is continuously differentiable}$

$\quad\quad\quad\quad\mbox{in}\,\, (\theta,\eta)\,\,\mbox{for}\,\,
\lambda-\mbox{almost all}\,\, y$

$\quad$

$\langle \mathbf{A2} \rangle \quad\quad \bigr|\dot{l}\,(\theta,
\eta) \bigr| \in L_2 (\R, G_{\theta, \eta})\,\,\,\mbox{for
all}\,\, (\theta, \eta) \in \Theta \times \Lambda $

$\quad$

$\langle \mathbf{A3} \rangle \quad\quad I(\theta, \eta)
\,\,\mbox{is nonsingular for all}\,\, (\theta, \eta) \in \Theta
\times \Lambda \,\,\mbox{and continuous}$

$\quad\quad\quad\quad\mbox{in}\, (\theta,\eta)\,.$
\end{definition}


\noindent This is a joint regularity condition and it is stronger
than the assumption that the model is regular in $\theta$ and
$\eta$ separately. Let us write $I(\theta_0, \eta_0)$ in the block
matrix form:

\begin{equation}\label{CompositeDeconvolution6}
I(\theta_0, \eta_0)\,= \,\left(%
\begin{array}{cc}
  I_{11}(\theta_0, \eta_0) & I_{12}(\theta_0, \eta_0) \\
  I_{21}(\theta_0, \eta_0) & I_{22}(\theta_0, \eta_0) \\
\end{array}%
\right)\,,
\end{equation}

\noindent where $I_{11}(\theta_0, \eta_0)$ is $k \times k,$
$I_{12}(\theta_0, \eta_0)$ is $k \times m,$ $I_{21}(\theta_0,
\eta_0)$ is $m \times k,$ $I_{11}(\theta_0, \eta_0)$ is $m \times
m.$ Thus, denoting for simplicity of formulas $\Omega:=
[y:\,g\,(y\,;(\theta_0,\eta_0))>0]$ we can write explicitly

\begin{equation}\label{CompositeDeconvolution7}
I_{11}(\theta_0, \eta_0)\,=\,E_{\theta_0, \eta_0}
\dot{l}_{\theta_0}^T\,\dot{l}_{\theta_0}\,=\,\int_{\R}
{\dot{l}_{\theta_0}}^T(y)\,\dot{l}_{\theta_0}(y)\, \,dG_{\theta_0,
\eta_0}(y)\,
\end{equation}

$$
=\,\int_{\Omega} \frac{
{\frac{\,\partial}{\,\partial\theta}\,\Bigr(
\int_{\R}\,f_{\theta}(s)\,h_{\eta_0}(\,y-s)\,ds
\Bigr)}^T\Bigr|_{\theta=\theta_0}
\frac{\,\partial}{\,\partial\theta}\,\Bigr(
\int_{\R}\,f_{\theta}(s)\,h_{\eta_0}(\,y-s)\,ds
\Bigr)\Bigr|_{\theta=\theta_0}
}{\int_{\R}\,f_{\theta_0}(s)\,h_{\eta_0}(\,y-s)\,ds}\,dy\,,
$$

\begin{equation}\label{CompositeDeconvolution8}
I_{12}(\theta_0, \eta_0)\,=\,E_{\theta_0, \eta_0}
\dot{l}_{\theta_0}^T\,\dot{l}_{\eta_0}\,=\,\int_{\R}
{\dot{l}_{\theta_0}}^T(y)\,\dot{l}_{\eta_0}(y)\, \,dG_{\theta_0,
\eta_0}(y)\,
\end{equation}

$$
=\,\int_{\Omega} \frac{ {\frac{\,\partial
}{\,\partial\theta}\,\Bigr(
\int_{\R}\,f_{\theta}(s)\,h_{\eta_0}(\,y-s)\,ds
\Bigr)}^T\Bigr|_{\theta=\theta_0}
\frac{\,\partial}{\,\partial\eta}\,\Bigr(
\int_{\R}\,f_{\theta_0}(s)\,h_{\eta}(\,y-s)\,ds
\Bigr)\Bigr|_{\eta=\eta_0}
}{\int_{\R}\,f_{\theta_0}(s)\,h_{\eta_0}(\,y-s)\,ds}\,dy\,,
$$



\noindent and analogously for $ I_{21}(\theta_0, \eta_0)$ and
$I_{22}(\theta_0, \eta_0).$ The \emph{efficient score function}
for $\theta$ in $M_{k,m}(\theta, \eta)$ is defined as (see
\cite{MR1245941}, p.28):

\begin{equation}\label{CompositeDeconvolution11}
l^{\ast}_{\theta_0}(y)\,=\, \dot{l}_{\theta_0}(y)\,-\,
I_{12}(\theta_0, \eta_0)\, I_{22}^{-1}(\theta_0,
\eta_0)\,\dot{l}_{\eta_0}(y)\,,
\end{equation}

\noindent and the \emph{efficient Fisher information matrix} for
$\theta$ in $M_{k,m}(\theta, \eta)$ is defined as

\begin{equation}\label{CompositeDeconvolution12}
I^{\ast}_{\theta_0}\,=\,E_{\theta_0, \eta_0}
l^{\ast\,T}_{\theta_0}\,l^{\ast}_{\theta_0}\,=\,\int_{\R}
l^{\ast}_{\theta_0}(y)^T\,l^{\ast}_{\theta_0}(y) \,dG_{\theta_0,
\eta_0}(y)\,.
\end{equation}

\noindent Before closing this section we consider two simple
examples.

\noindent {\bf Example 2.} Suppose $\theta \in \R,$ $\eta \in
\R^+$ and, moreover, $\{f_\theta\}$ is a family $\{N(\theta, 1)\}$
of normal densities with mean $\theta$ and variance 1, and
$\{h_\eta\}$ is a family $\{N(0, \eta^2 )\}.$ Then
$g(\theta,\eta)= f_\theta \ast h_\eta \sim N(\theta, \eta^2 +1).$
Let $\theta$ be the parameter of interest and $\eta$ the nuisance
one. Let $H_0$ be $\theta=\theta_0.$ By
(\ref{CompositeDeconvolution3}) and
(\ref{CompositeDeconvolution4}) for all $y$

\begin{equation}\label{CompositeDeconvolution19}
\dot{l}_{\theta_0}(y)= \, \frac{y-\theta_0}{\eta_0^2 +1}\,, \quad
\dot{l}_{\eta_0}(y)= \,\frac{ (y-\theta_0)^2 \,\eta_0}{{(\eta_0^2
+1)}^{2}} - \frac{\eta_0}{\eta_0^2 +1} \,.
\end{equation}








\noindent By (\ref{CompositeDeconvolution8})

$$I_{12}(\theta, \eta)\,=\,\int_{\R} \frac{y-\theta}{\eta^2 +1} \biggr[ \frac{ (y-\theta)^2
\eta}{{(\eta^2 +1)}^{2}} - \frac{\eta}{\eta^2 +1}
\biggr]\,dN(\theta, \eta^2 +1)(y)\,=\,0\,,$$




\noindent for all $\theta, \eta.$ This means that adaptive
estimation of $\theta$ is possible in this model, i.e., we can
estimate $\theta$ equally well whether we know the true $\eta_0$
or not. Though, we will not be concerned with estimation here.
From (\ref{CompositeDeconvolution7}) we get

\begin{equation}\label{CompositeDeconvolution21}
{(I^\ast_\theta)}^{-1}\,=\,\int_{\R} \frac{(y-\theta)^2}{{(\eta^2
+1)}^{2}} \,dN(\theta, \eta^2 +1)(y)\,=\, \frac{1}{\eta^2 +1} \,.
\end{equation}


\noindent {\bf Example 3.} Suppose now that we are interested in
the parameter $\eta$ in the situation of Example 2 and the null
hypothesis is $H_0: \eta=\eta_0.$ There is a sort of symmetry
between signal and noise: "what is a signal for one person is a
noise for the other" (see also Remark
\ref{CompositeDeconvolutionRemark1}). From Example 2 we know that
the score function $\dot{l}_{\eta_0}$ for $\eta$ at $\eta_0$ is
given by (\ref{CompositeDeconvolution19}). Since we proved for
this example $I_{12}=I_{21}=0,$ the efficient score function
$l^{\ast}_{\eta_0}$  for $\eta$ at $\eta_0$ is given by
(\ref{CompositeDeconvolution19}) as well. We calculate now

\begin{equation}\label{CompositeDeconvolution24}
{(I^\ast_{\eta_0})}^{-1}\,=\,\int_{\R} {\biggr( \frac{
(y-\theta)^2 \,\eta_0}{{(\eta_0^2 +1)}^{2}} -
\frac{\eta_0}{\eta_0^2 +1} \biggr)}^2 \,dN(\theta, \eta_0^2
+1)(y)\,=:\,\frac{1}{C(\eta_0)}\,.
\end{equation}



\noindent The constant $C(\eta_0)$ in
(\ref{CompositeDeconvolution24}) can be expressed explicitly in
terms of $\eta_0,$ but this is not the point of this example. By
the symmetry of $\theta$ and $\eta$ we have
$l^{\ast}_{\eta_0}(y)\,=\, \dot{l}_{\eta_0}(y)\,-\, I_{21}(\theta,
\eta_0)\, I_{11}^{-1}(\theta,
\eta_0)\,\dot{l}_{\theta_0}(y)\,=\,\dot{l}_{\eta_0}(y)\,.$

\begin{remark}\label{CompositeDeconvolutionRemark1} Note that the
problem is symmetric in $\theta$ and $\eta$ in the sense that it
is possible to consider estimating and testing for each parameter,
$\theta$ or $\eta.$ Physically this means that from the noisy
signal one can recover some "information" not only about the pure
signal but also about the noise. This is actually natural since a
noise is in fact also a signal. We are observing two signals at
once. The payment for this possibility is that except for some
trivial cases one can't recover full information about both the
signal of interest as well as about the noise.
\end{remark}

\section{Efficient score test}\label{CompositeDeconvolutionSection4}

Let $l^{\ast}_{\theta_0}$ be defined by
(\ref{CompositeDeconvolution11}) and $I^{\ast}_{\theta_0}$ by
(\ref{CompositeDeconvolution12}). Note that both
$l^{\ast}_{\theta_0}$ and $I^{\ast}_{\theta_0}$ depends (at least
in principle) on the unknown nuisance parameter $\eta_0.$ Let
$l^{\ast}_j$ and $L$ be some estimators of
$l^{\ast}_{\theta_0}(Y_j)$ and ${(I^{\ast}_{\theta_0})}^{-1}$
correspondingly. These estimators are supposed to depend only on
the observable $Y_1, \ldots , Y_n\,,$ but not on the $X_1, \ldots
, X_n.$

\begin{definition}\label{CompositeDeconvolutionDef2}
We say that $l_j^\ast$ is a \emph{sufficiently good} estimator of
$l^{\ast}_{\theta_0}(Y_j)$ if for each $(\theta_0, \eta_0) \in
\Theta \times \Lambda$ it holds that for every $\varepsilon
>0$
\begin{equation}\label{CompositeDeconvolution14}
G^{n}_{\theta_0, \eta_0} \biggr(\frac{1}{\sqrt{n}} \,\biggr\|
\sum_{i=1}^{n}(l_j^\ast - l^{\ast}_{\theta_0}(Y_j)) \biggr\| \geq
\varepsilon \biggr ) \rightarrow 0 \quad\mbox{as}\quad n
\rightarrow\infty\,,
\end{equation}

\noindent where $\| \cdot \|$ denotes the Euclidian norm of a
given vector.
\end{definition}

\noindent In other words, condition
(\ref{CompositeDeconvolution14}) means that the average
$\frac{1}{n}\,\sum_{i=1}^{\,n}l^{\ast}_{\theta_0}(Y_j)\,\approx \,
E_{\theta_0, \eta_0} l^{\ast}_{\theta_0}$ is
$\sqrt{\,n}-$consistently estimated. We illustrate this definition
by some examples.

\noindent{\bf Example 2 (continued)}. We have (denoting variance
of $Y$ by $\sigma^2 (Y)$):

$$l^{\ast}_{\theta_0}(Y_j) = \frac{Y_j - \theta_0}{\sigma^2 (Y)}.$$

\noindent Define $$l^{\ast}_j := \frac{Y_j -
\theta_0}{\widehat{\sigma}_n^2},$$ where $\widehat{\sigma}_n^2$ is
any $\sqrt{n}-$consistent estimator of the variance of $Y.$ One
can take, for example, the sample variance $s_n^2 = s_n^2 (Y_1,
\ldots , Y_n)$ as such an estimate. Then, since by the model
assumptions $\sigma^2 (Y) >0,$ thus constructed $l^{\ast}_j$
satisfies Definition \ref{CompositeDeconvolutionDef2}. See
Appendix for the proof. $\Box$

\noindent{\bf Example 3 (continued)}. We have in this case

$$l^{\ast}_{\eta_0}(Y_j) = \frac{\eta_0}{\eta_0^2 +1}(Y_j - \theta_0)^2 - \frac{\eta_0}{\eta_0^2
+1}.$$ For simplicity of notations we write
$l^{\ast}_{\eta_0}(Y_j) = C_1(\eta_0)(Y_j - \theta_0)^2 -
C_2(\eta_0).$ Let $\widehat{\theta}_n$ be any
$\sqrt{n}-$consistent estimate of $\theta_0$ and put $l^{\ast}_j
:= C_1(\eta_0)(Y_j - \widehat{\theta}_n)^2 - C_2(\eta_0).$ Then
Definition \ref{CompositeDeconvolutionDef2} is satisfied in this
Example also. This is proved in Appendix. $\Box$

\noindent Definition \ref{CompositeDeconvolutionDef2} reflects the
basic idea of the method of estimated scores. This method is
widely used in statistics (see \cite{MR1245941}, \cite{MR856811},
\cite{MR620321}, \cite{MR2259976} and others). These authors show
that for different problems it is possible to construct nontrivial
parametric, semi- and nonparametric estimators of scores such that
these estimators will satisfy (\ref{CompositeDeconvolution14}).

\begin{definition}\label{CompositeDeconvolutionDef3}

Define

\begin{equation}\label{CompositeDeconvolution15}
W_k =\biggr\{ \frac{1}{\sqrt{n}}\, \sum_{j=1}^{n} l^\ast_j
\biggr\} \;\widehat{L}\; {\biggr\{ \frac{1}{\sqrt{n}}\,
\sum_{j=1}^{n} l^\ast_j \biggr\} }^{T} \,,
\end{equation}

\noindent where $\widehat{L}$ is an estimate of
${(I^{\ast}_{\theta_0})}^{-1}$ depending only on $Y_1, \ldots ,
Y_n.$ Note that $l^\ast_j$ is a $k-$dimensional vector and
$\widehat{L}$ is a $k \times k$ matrix. We call $W_k$ the
\emph{efficient score test statistic} for testing $H_0:
\theta=\theta_0$ in $M_{k,m}(\theta, \eta).$ It is assumed that
the null hypothesis is rejected for large values of $W_k.$
\end{definition}

Normally it should be possible to construct reasonably good
estimators $\widehat{\eta}_n$ of $\eta$ by standard methods since
at this point our construction is parametric. After that it would
be enough to plug in these estimates in
(\ref{CompositeDeconvolution11}) and get the desired
${l^\ast}_j'$s satisfying (\ref{CompositeDeconvolution14}).

\noindent {\bf Example 2 (continued).} Let $\widehat{\sigma}^2
(Y)$ be any $\sqrt{n}-$consistent estimate of $\eta^2 +1$ such
that this estimate is based on $Y_1, \ldots , Y_n.$ Then by
(\ref{CompositeDeconvolution21}), (\ref{CompositeDeconvolution19})
and definition (\ref{CompositeDeconvolution15}) the efficient
score test statistic for testing $H_0: \theta=\theta_0$ (in the
model $M_{1,1}(\theta, \eta)$) is

\begin{equation}\label{CompositeDeconvolution22}
W_1 =\biggr( \frac{1}{\sqrt{n}}\, \sum_{j=1}^{n}
\frac{Y_j-\theta_0}{ \widehat{\sigma}^2_n (Y) } \biggr)^2 \,
\widehat{\sigma}_n^2 (Y)\,=\, \frac{1}{\widehat{\sigma}^2_n (Y)}
 \biggr( \frac{1}{\sqrt{n}}\,
\sum_{j=1}^{n} (Y_j-\theta_0)\biggr)^2 \,.
\end{equation}

\noindent {\bf Example 3 (continued).} Using any $\sqrt{n}-$
consistent estimate $\widehat{\theta}$ of $\theta,$ we get the
efficient score test statistic

\begin{eqnarray}
 W_1 & = & \biggr( \frac{1}{\sqrt{n}}\, \sum_{j=1}^{n}
 \,\biggr[\frac{ (Y_j-\widehat{\theta}_n)^2 \,\eta_0}{{(\eta_0^2
+1)}^{2}} - \frac{\eta_0}{\eta_0^2 +1}\biggr]
\,\biggr)^2\,C(\eta_0) \nonumber\\
     & = & \biggr( \frac{1}{\sqrt{n}}\, \frac{\eta_0}{{(\eta_0^2
+1)}^{2}} \,\sum_{j=1}^{n} \,(Y_j-\widehat{\theta}_n)^2 \,-\,
\sqrt{n}\,\frac{\eta_0}{\eta_0^2 +1}\,\biggr)^2\,C(\eta_0) \,.
\end{eqnarray}



\begin{remark}\label{CompositeDeconvolutionRemark3} We make now
the following remark to avoid possible confusions. For the simple
deconvolution we had the score test statistics and now we have the
\emph{efficient} score test statistics. This does not mean that
the statistics for simple deconvolution is "inefficient". Here the
word "efficient" has a strictly technical meaning. Because of the
presence of the nuisance parameter we have to extract information
about the parameter of interest. We want to do this efficiently in
some sense. This is the explanation of the terminology.
\end{remark}

The following theorem describes asymptotic behavior of $W_k$ under
the null hypothesis.

\begin{theorem}\label{CompositeDeconvolutionTh1}
Assume the null hypothesis $H_0: \theta=\theta_0$ holds true,
$\langle A1 \rangle$-$\langle A3 \rangle$ are fulfilled,
(\ref{CompositeDeconvolution14}) is satisfied, and $\widehat{L}$
is any consistent estimate of ${(I^{\ast}_{\theta_0})}^{-1}.$ Then

$$W_k \rightarrow_d \,{\chi}^{2}_{k} \quad\mbox{as}\quad n \rightarrow \infty\,,$$ where
${\chi}^{2}_{k}$ denotes a random variable with central chi-square
distribution with $k$ degrees of freedom.
\end{theorem}
\vspace{6pt}

\section{Selection rule}\label{CompositeDeconvolutionSection5}


In this section we extend the construction of Section
\ref{SimpleDeconvolutionSubsection3.1} to the case of composite
hypotheses. First we give a general definition of a selection
rule.

Denote by $M_{k,m} (\theta, \eta)$ the model described in Section
\ref{CompositeDeconvolutionSection2} and such that the true
parameter $(\theta, \eta)$ belongs to a parameter set, say
$\Theta_k \times \Lambda ,$ and $\dim \Theta_k = k.$ By a
\emph{nested family} of submodels $M_{k,m} (\theta, \eta)$ for
$k=1, \ldots$ we would mean a sequence of these models such that
for their parameter sets it holds that $\Theta_1 \times \Lambda
\subset \Theta_2 \times \Lambda \subset \ldots .$

\begin{definition}\label{CompositeDeconvolutionDef4}
Consider a nested family of submodels $M_{k,m}(\theta, \eta)$ for
$k=1,$ $\ldots ,$ $d,$ where $d$ is fixed but otherwise arbitrary,
and $m$ is fixed. Choose a function $\pi (\cdot,\cdot): \mathbb{N}
\times \mathbb{N} \rightarrow \mathbb{R},$ where $\mathbb{N}$ is
the set of natural numbers. Assume that $\pi(1,n) <\pi(2,n)<
\ldots <\pi(d,n)$ for all $n$ and $\pi(j,n)-\pi(1,n) \rightarrow
\infty$ as $n \rightarrow \infty$ for every $j= 2, \ldots ,d.$
Call $\pi(j,n)$ a \emph{penalty attributed to the j-th model
$M_j(\theta)$ and the sample size n.} Then a \emph{selection rule}
$S(l^\ast)$ for the test statistic $W_k$ is an integer-valued
random variable satisfying the condition

\begin{equation}\label{CompositeDeconvolution17}
S(l^\ast)=\min\bigr\{ k:\,1\leq k \leq d;\, W_k - \pi(k,n) \geq
W_j -\pi(j,n),\, j= 1, \ldots ,d\bigr\}\,.
\end{equation}

\noindent We call the random variable $W_S$ a \emph{data-driven
efficient score test statistic} for testing validity of the
initial model. In this paper we also assume that the following
condition holds.

$\quad$

$\langle \mathbf{S1} \rangle \quad\quad\mbox{for every fixed }\,\,
k \geq 1 \,\,\mbox{it holds that}\,\, \pi (k,n)= o(n)
\,\,\mbox{as}\,\, n \rightarrow \infty\,.$

$\quad$
\end{definition}

Unlike the case of the simple null hypothesis, in the case of the
composite hypotheses the selection rule depends on the estimator
$l^\ast_j$ of the unknown values $l^{\ast}_{\theta_0}(Y_j)$ of the
efficient score function. This means that we need to estimate the
nuisance parameter $\eta ,$ or corresponding scores, or their sum.
Surprising result follows from Theorem
\ref{CompositeDeconvolutionTh2} below: for our problem many
possible penalties and, moreover, essentially all sensible
estimators plugged in $W_k,$ give consistent selection rules.
Possible choices of penalties are, for instance, Shwarz's penalty
$\pi(j,n)=j \, \log n,$ or Akaike's penalty $\pi(j,n)=j.$

Denote by $P^n_{\theta_0, \eta_0}$ the probability measure
corresponding to the case when $X_1, \ldots ,$ $X_n$ all have the
density $f(\theta_0, \eta_0).$ The main result about the
asymptotic null distribution of $W_S$ is the following theorem (it
is proved analogously to Theorem \ref{SimpleDeconvolutionTh2}).

\begin{theorem}\label{CompositeDeconvolutionTh2}
Under the conditions of Theorem \ref{CompositeDeconvolutionTh1},
as $n \rightarrow \infty$ it holds that
$$P^n_{\theta_0, \eta_0} (S(l^\ast) >1) \rightarrow 0 \quad \mbox{and}\quad W_S
\rightarrow_d \,{\chi}^{2}_{1}.$$
\end{theorem}


\noindent Condition (\ref{CompositeDeconvolution14}) is what makes
this direct reference to the case of the simple hypothesis
possible. Estimation of the efficient score function
$l^{\ast}_{\theta_0}$ can be done by different ways. First way is
to estimate the whole expression from the right side of
(\ref{CompositeDeconvolution11}). For this method of estimation
condition (\ref{CompositeDeconvolution14}) is natural. The second
and probably more convenient method of estimating
$l^{\ast}_{\theta_0}$ is via estimation of the nuisance parameter
$\eta$ by some estimator $\widehat{\eta}.$ But for this approach
condition (\ref{CompositeDeconvolution14}) becomes something that
have to be proved for each particular estimator. We hope that this
inconvenience is excused by the fact that we are only introducing
the new test here. It is possible to reformulate condition
(\ref{CompositeDeconvolution14}) explicitly in terms of conditions
on $\widehat{\eta},$ $\{f_\theta\},$ and $\{h_\eta\}$ (see an
analogue in \cite{MR1447749}).

\begin{remark}\label{CompositeDeconvolutionRemark4} The selection rule
$S(l^\ast)$ can be modified in order to make it possible to choose
not only models of dimension less than some fixed $d,$ but to
allow arbitrary large dimensions of $M_{k,m}(\theta, \eta)$ as the
number of observations grows. See Remark
\ref{SimpleDeconvolutionRemark1}.
\end{remark}

\begin{remark}\label{CompositeDeconvolutionRemark5} It is possible to
modify the definition of selection rule so that both dimensions
$k$ and $m$ would be selected by the test from the data. A
corresponding test statistic will be of the form $W_S,$ where this
time $S=(S_1, S_2).$ Proofs of the asymptotic properties for this
statistic are analogous to those presented in this paper. Possibly
this statistic could be useful since the situation with the noise
of an unknown dimension often seems to be more realistic. On the
other hand, this statistic will also have some disadvantages. One
will have to impose more strict assumptions on both signal and
noise (including an analogue of the double-identifiability
assumption). Also the final result will be weaker than the result
of this section. This will be a payment for an attempt to extract
information about a larger number of parameters from the same
amount of observations $Y_1,$ $\ldots ,$ $Y_n\,.$
\end{remark}

\section{Consistency of tests}\label{CompositeDeconvolutionSection6}

Let $F$ be a true distribution function of $X$ and $H$ a true
distribution of $\varepsilon.$ Here $F$ and $H$ are \emph{not}
necessarily parametric and possibly these distribution functions
do not have densities with respect to the Lebesgue measure
$\lambda .$ Let us choose for every $k \leq d$ an auxiliary
parametric family $\{f_{\theta}\},\,$ $\theta \in \Theta \subseteq
\mathbb{R}^k$ such that $f_0$ from this family coincides with
$f_0$ from the null hypothesis $H_0.$ Correspondingly, let us fix
an integer $m$ and choose an auxiliary parametric family
$\{h_{\eta}\},\,$ $\eta \in \Lambda \subseteq \mathbb{R}^m.$
Suppose that the chosen families $\{f_{\theta}\}$ and
$\{h_{\eta}\}$ give us the regular deconvolution problem in the
sense of Definition \ref{CompositeDeconvolutionDefinition1}. Then
one is able to construct the score test statistic $W_k$ defined by
(\ref{CompositeDeconvolution15}) despite the fact that the true
$F$ and $H$ possibly do not have any relation to the chosen
$\{f_{\theta}\}$ and $\{h_{\eta}\}.$ This is our goal in this
section to determine under what conditions thus build $W_k$ will
be consistent for testing against $H_A.$

Suppose that the following condition holds

$\quad$

$\langle \mathbf{C1} \rangle \quad\quad\mbox{there exists
integer}\,\, K \geq 1 \,\,\mbox{such that}\,\, K \leq d
\,\,\mbox{and}$

$\quad\quad\quad\quad E_{F \ast H}\,l^{\ast}_{\theta_0(1)}=0,
\ldots , E_{F \ast H}\,l^{\ast}_{\theta_0(K-1)}=0,\,\,E_{F \ast
H}\,l^{\ast}_{\theta_0(K)}= C_K \neq 0\,,$

$\quad$

\noindent where $l^{\ast}_{\theta_0(i)}$ is the $i-$th coordinate
function of $l^{\ast}_{\theta_0}$ and $l^{\ast}_{\theta_0}$ is
defined by (\ref{CompositeDeconvolution11}), $d$ is the maximal
possible dimension of our model as in Definition
\ref{CompositeDeconvolutionDefinition1} of Section
\ref{CompositeDeconvolutionSection5}, and $E_{F \ast H}$ denotes
the mathematical expectation with respect to $F \ast H.$

Condition $\langle C1 \rangle$ is a weak analog of nondegeneracy:
if for all $k$ $\langle C1 \rangle$ fails, then $F $ is orthogonal
to the whole system ${l^{\ast}_{\theta_0(i)}}_{i=1}^{\infty}$ and
if this system is complete, then $F \ast H$ is degenerate. Also
$\langle C1 \rangle$ is related to the identifiability of the
model: if the model is not identifiable, then $F \ast H = F_0 \ast
H$ can happen and $\langle C1 \rangle$ fails. Establishing
identifiability for the parametric deconvolution is not trivial
(see \cite{MR0270486}, e.g.). It is important to note also that
although $\langle C1 \rangle$ has something common with both
nondegeneracy and identifiability, it is in general pretty far
from both these notions.


The main result of this section is the following.

\begin{theorem}\label{CompositeDeconvolutionTh3}
If (\ref{CompositeDeconvolution14}) is satisfied and $\widehat{L}$
is any consistent estimate of ${(I^{\ast}_{\theta_0})}^{-1},$ then

\begin{enumerate}
\item \noindent the test based on $W_k$ is consistent for testing
against all alternative distributions $F,$ $H$ such that $\langle
C1 \rangle$ is satisfied with $K \leq k$

\item \noindent the test based on $W_k$ is inconsistent for
testing against alternative distributions $F,$ $H$ such that
$\langle C1 \rangle$ is satisfied with $K > k$

\item \noindent if the selection rule $S(l^\ast)$ satisfies
$\langle S1 \rangle,$ then test based on $W_S$ is consistent
against all alternative distributions $F \ast H$ such that
$\langle C1 \rangle$ is satisfied with some $K.$
\end{enumerate}
\end{theorem}

\noindent Part 2 of Theorem \ref{CompositeDeconvolutionTh3} shows
why it is important to choose the suitable model dimension. Now we
give two specific examples.

\noindent {\bf Example 2 (continued).} By Theorem
\ref{CompositeDeconvolutionTh3} the test based on $W_1$ is
consistent if and only if for true $F$ and $H$ it holds that

\begin{equation}\label{CompositeDeconvolution23}
\frac{1}{\eta^2 +1}\,E_{F \ast H}(Y - \theta_0) \neq 0\,, \quad
\mbox{i.e.} \quad E_{F \ast H}(Y) \neq \theta_0\,.
\end{equation}

\noindent For example, $W_1$ doesn't work when the true $H$ is
symmetric about 0 and the true $F \neq F_0$ has the mean equal to
$\theta_0.$

\noindent {\bf Example 3 (continued).} By Theorem
\ref{CompositeDeconvolutionTh3} $W_1$ is consistent if and only if
for true $F$ and $H$ it holds that

$$ E_{F \ast H} \,\biggr[\frac{ (y-\theta)^2 \,\eta_0}{{(\eta_0^2
+1)}^{2}} - \frac{\eta_0}{\eta_0^2 +1} \biggr]\, \neq 0\,, \quad
\mbox{i.e.}$$

\begin{equation}\label{CompositeDeconvolution26}
E_{F \ast H} \,(y-\theta)^2 \, \neq \eta_0^2 +1, \quad \mbox{or
equivalently} \quad Var_{F \ast H}\,Y \,\neq\, Var_{F \ast
H_0}\,Y\,.
\end{equation}


\noindent Note that condition (\ref{CompositeDeconvolution23}) can
be interpreted as "$W_1$ is consistent for testing the hypothesis
about the mean in this model iff the expectation of $Y$ under
alternative is different from the expectation under the null
hypothesis" and (\ref{CompositeDeconvolution26}) as "$W_1$ is
consistent for testing the hypothesis about the variance in this
model iff the variance of $Y$ under alternative is different from
the variance under the null hypothesis". One cannot expect more
from such a simple test as $W_1$. On contrary, the data-driven
test statistic $W_S$
provides a consistent testing procedure.\\

\section*{Acknowledgments} Author would like to thank Axel Munk
for suggesting this topic of research and Fadoua Balabdaoui,
Ta-Chao Kao, Wilbert Kallenberg and Janis Valeinis for
helpful discussions. 

\bibliographystyle{plainnat}
\bibliography{Deconvolution_Bibliography}

\noindent {\bf Appendix.}

\begin{proof}(Theorem \ref{SimpleDeconvolutionTh3}).
We calculated the efficient score vector in
(\ref{SimpleDeconvolution16})-(\ref{SimpleDeconvolution17}). By
Proposition 1, p.13 of \cite{MR1245941} and our regularity
assumptions matrix $L$ exists and is positive definite and
nondegenerate of rank $k.$ Under $\langle B1 \rangle - \langle B3
\rangle$ $E_0 l^\ast (y)=0$ (see \cite{MR1245941}, p.15) and our
statement follows.
\end{proof}

\begin{proof}(Proposition \ref{SimpleDeconvolutionProp5}).
Follows by the multivariate Central Limit Theorem.
\end{proof}

\begin{proof}(Theorem \ref{SimpleDeconvolutionTh2}).
Denote $\Delta (k,n) := \pi(k,n)-\pi(1,n).$ For any $k = 2, \ldots
,d$

\begin{eqnarray}
 P^n_0 (S=k) & \leq & P^n_0 \bigr( U_k - \pi(k,n) \geq U_1
-\pi(1,n)\bigr) \nonumber\\
             & \leq & P^n_0 \bigr( U_k \geq
\pi(k,n)-\pi(1,n)\bigr) \nonumber\\
             & = & P^n_0 \bigr( U_k \geq \Delta(k,n)\bigr)\,.\nonumber
\end{eqnarray}



\noindent By Theorem \ref{SimpleDeconvolutionTh3} $U_k
\rightarrow_d \,{\chi}^{2}_{k}$ as $n \rightarrow \infty,$ thus
for $\Delta (k,n) \uparrow \infty$ as $n \rightarrow \infty$ we
have $P^n_0 \bigr(U_k \geq \Delta (k,n)\bigr) \rightarrow 0$ as $n
\rightarrow \infty,$ so for any $k = 2, \ldots ,d$ we have $P^n_0
(S=k) \rightarrow 0$ as $n \rightarrow \infty.$ This proves that

$$P^n_0 (S \geq 2) = \sum_{k=2}^{d} P^n_0 (S=k) \rightarrow 0, \quad n \rightarrow \infty
,$$ and so $P^n_0 (S=1) \rightarrow 1.$ Now write for arbitrary
real $t > 0$

\begin{eqnarray}\label{SimpleDeconvolution12}
 P^n_0 (|U_S - U_1| \geq t) & = & P^n_0 (|U_1 - U_1| \geq t; \, S=1) \nonumber\\
                            &   & \quad + \sum_{m=2}^{d} P^n_0 (|U_m - U_1| \geq t; \,
             S=m)\nonumber\\
                            & = & \sum_{m=2}^{d} P^n_0 (|U_m - U_1| \geq t; \, S=m).
\end{eqnarray}



\noindent For $m = 2, \ldots ,d$ we have $P^n_0 (S=m) \rightarrow
0,$ so

$$0 \leq \sum_{m=2}^{d} P^n_0 (|U_m - U_1| \geq t; \, S=m) \leq \sum_{m=2}^{d}
P^n_0 (S=m) \rightarrow 0$$ as $n \rightarrow \infty$ and thus by
(\ref{SimpleDeconvolution12}) it follows that $U_S$ tends to $U_1$
in probability as $n \rightarrow \infty.$ But $U_1 \rightarrow_d
\,{\chi}^{2}_{1}$ by Theorem \ref{SimpleDeconvolutionTh3}, so $U_S
\rightarrow_d \,{\chi}^{2}_{1}$ as $n \rightarrow \infty.$
\end{proof}

We shall use the following standard lemma from linear algebra.

\begin{lemma}\label{SimpleDeconvolutionLemmaAlgebra1}
Let $x \in {\R}^k,$ and let $A$ be a $k \times k$ positive
definite matrix; if for some real number $\delta >0$ we have $A >
\delta$ (in the sense that the matrix $(A - \delta\, I_{k \times
k})$ is positive definite, where $I_{k \times k}$ is the $k \times
k$ identity matrix), then for all $x \in {\R}^k$ it holds that $xA
x^T > \delta {\|x\|}^2.$
\end{lemma}

\begin{proof}(Proposition \ref{SimpleDeconvolutionProp6}).
From $\langle D1 \rangle$ by the law of large numbers we get

\begin{equation}\label{SimpleDeconvolution21}
\frac{1}{n}\,\sum_{j=1}^{n} l_i^\ast (Y_j)\,\rightarrow_P \,0
\quad\mbox{for}\quad 1 \leq i \leq K-1
\end{equation}

\begin{equation}\label{SimpleDeconvolution22}
\frac{1}{n}\,\sum_{j=1}^{n} l_i^\ast (Y_j)\,\rightarrow_P \,C_K
\neq 0.
\end{equation}

\noindent We apply Lemma \ref{SimpleDeconvolutionLemmaAlgebra1} to
the matrix $L$ defined in (\ref{SimpleDeconvolution18}); since all
the eigenvalues of $L$ are positive we can choose $\delta$ to be
any fixed positive number less than the smallest eigenvalue of
$L.$ We obtain the following inequality

\begin{eqnarray}\label{SimpleDeconvolution23}
 U_k & = & \biggr\{ \frac{1}{\sqrt{n}}\, \sum_{j=1}^{n} l^\ast (Y_j)
\biggr\} \;L\; {\biggr\{ \frac{1}{\sqrt{n}}\, \sum_{j=1}^{n}
l^\ast (Y_j) \biggr\} }^{T} \nonumber\\
     & > & \delta {\biggr\|
\frac{1}{\sqrt{n}}\, \sum_{j=1}^{n} l^\ast (Y_j) \biggr\|}^2 =
\delta\,n\, \sum_{i=1}^{k} {\biggr(\frac{1}{n}\,\sum_{j=1}^{n} l_i^\ast (Y_j) \biggr)}^{\!\!2}\nonumber\\
     & \geq & \delta\,n {\biggr(\frac{1}{n}\,\sum_{j=1}^{n} l_K^\ast
(Y_j) \biggr)}^{\!\!2} \,.
\end{eqnarray}



\noindent Now by (\ref{SimpleDeconvolution21}) and
(\ref{SimpleDeconvolution22}) we get for all $s \in \R$

\begin{eqnarray}
 P(U_k \leq s) & \leq & P \biggr(\delta\,n
{\biggr(\frac{1}{n}\,\sum_{j=1}^{n} l_K^\ast (Y_j)
\biggr)}^{\!\!2} \leq
s \biggr) \nonumber\\
               &   =  & P \biggr( {\biggr(\frac{1}{n}\,\sum_{j=1}^{n} l_K^\ast (Y_j)
\biggr)}^{\!\!2} \leq \frac{s}{\,\delta\,n} \biggr)\nonumber\\
               &   =  & P \biggr(\biggr|\frac{1}{n}\,\sum_{j=1}^{n} l_K^\ast (Y_j)
\biggr| \leq \sqrt{\frac{s}{\,\delta\,n}}\,\, \biggr) \rightarrow
0 \quad\mbox{as}\quad n \rightarrow \infty\,,\nonumber
\end{eqnarray}

\noindent and this proves the Proposition.


\end{proof}

\begin{proof}(Proposition \ref{SimpleDeconvolutionProp7}).
Let $\pi (k,n)$ and $\Delta (k,n)$ be defined as in Section
\ref{SimpleDeconvolutionSubsection3.1}. For any $i=1, \ldots ,
K-1$ we have

\begin{eqnarray}\label{SimpleDeconvolution24}
 P_F\,(S=i) & \leq & P_F\,\bigr(U_i - \pi (i,n) \geq U_K - \pi
(K,n)\bigr)\nonumber\\
            &   =  & P_F\,\bigr(U_i \geq U_K - (\pi (K,n) - \pi
(i,n))\bigr).
\end{eqnarray}



\noindent By (\ref{SimpleDeconvolution22}) and
(\ref{SimpleDeconvolution23}) we get

\begin{equation}\label{SimpleDeconvolution25}
P_F\,\biggr(U_K \geq \delta \,\frac{C_K}{2}\,n \biggr)\rightarrow
\,1 \quad\mbox{as}\quad n \rightarrow \infty\,.
\end{equation}

\noindent Note that


\begin{equation}\label{SimpleDeconvolution26}
P_F\,\bigr(U_i \geq U_K - (\pi (K,n) - \pi (i,n))\bigr)
\quad\quad\quad\quad
\end{equation}

\vskip -0.7cm

\begin{eqnarray}
 \quad\quad & \leq & P_F \biggr(U_i \geq \delta \,\frac{C_K}{2}\,n - (\pi (K,n) -
\pi (i,n));\,U_K \geq \delta \,\frac{C_K}{2}\,n \biggr)\nonumber\\
 \quad\quad & \quad & +\,\, P_F \biggr(U_K \leq \delta \,\frac{C_K}{2}\,n
 \biggr).\nonumber
\end{eqnarray}


\noindent Since by $\langle S1 \rangle$ it holds that $\pi (K,n) -
\pi (i,n)= o(n),$ we get

\begin{equation}\label{SimpleDeconvolution27}
P_F\,\biggr(U_i \geq \delta \,\frac{C_K}{2}\,n - (\pi (K,n) - \pi
(i,n));\,U_K \geq \delta \,\frac{C_K}{2}\,n \biggr)
\quad\quad\quad\quad\quad\quad\quad
\end{equation}

\vskip -0.5cm

\begin{eqnarray}
 \quad\quad & \leq & P_F\,\biggr(U_i \geq \delta \,\frac{C_K}{2}\,n - (\pi (K,n)
- \pi (i,n)) \biggr) \nonumber\\
 \quad\quad & \leq & P_F\,\biggr(U_i \geq \delta \,\frac{C_K}{2}\,n
\biggr)\rightarrow \,0\nonumber
\end{eqnarray}



\noindent as $n \rightarrow \infty$ by Chebyshev's inequality
since by Proposition \ref{SimpleDeconvolutionProp5} we have $U_i
\rightarrow_d \,{\chi}^{2}_{i} \quad\mbox{as}\quad n \rightarrow
\infty\,$ for all $i=1, \ldots , K-1.$ Substituting
(\ref{SimpleDeconvolution25}) and (\ref{SimpleDeconvolution27}) to
(\ref{SimpleDeconvolution26}) we get $P_F\,(S=i) \rightarrow \,0$
as $n \rightarrow \infty$ for all $i=1, \ldots , K-1.$ This means
that $P_F \,(S \geq K) \rightarrow 1$ as $n \rightarrow \infty.$

Now write for $t \in \R$

$$P_F\,(U_S \leq t) = P_F\,(U_S \leq t; S \leq K-1) + P_F\,(U_S \leq t; S \geq K)=: R_1 + R_2.$$

\noindent But $R_1 \rightarrow \,0$ since $P_F\,(S=i) \rightarrow
\,0$ for $i=1, \ldots , K-1$ and $K \leq d < \infty.$ Since
$U_{l_1} \geq U_{l_2}$ for $l_1 \geq l_2,$ we get

$$R_2 \leq \sum_{l=K}^{d}\,P_F\,(U_K \leq t) \rightarrow \,0$$ as
$n \rightarrow \infty$ by Proposition
\ref{SimpleDeconvolutionProp6}. Thus $P_F \,(U_S \leq t)
\rightarrow 0$ as $n \rightarrow \infty$ for all $t \in \R.$
\end{proof}

\begin{proof}(Theorem \ref{SimpleDeconvolutionTh8}).
Part 1 follows from Theorem \ref{SimpleDeconvolutionTh3} and
Proposition \ref{SimpleDeconvolutionProp6}, part 2 from Theorem
\ref{SimpleDeconvolutionTh3} and Proposition
\ref{SimpleDeconvolutionProp5}, part 3 from Theorem
\ref{SimpleDeconvolutionTh2} and Proposition
\ref{SimpleDeconvolutionProp7}.
\end{proof}

\begin{proof} (The statement about $l^{\ast}_j$ from Example 2).
Indeed,

\begin{equation*}
 \begin{split}
 \frac{1}{\sqrt{n}} \,\biggr| \sum_{i=1}^{n}(l_j^\ast -
  l^{\ast}_{\theta_0}(Y_j)) \biggr| & = & \frac{1}{\sqrt{n}} \,\biggr|
  \sum_{i=1}^{n}( Y_j - \theta_0 )\Bigr( \frac{1}{\sigma^2 (Y)} -
  \frac{1}{\widehat{\sigma}_n^2} \Bigr) \biggr| \\
                                    & = & \sqrt{n} \,\biggr| \frac{1}{\sigma^2 (Y)} -
\frac{1}{\widehat{\sigma}_n^2} \biggr| \cdot \frac{1}{n} \,\biggr|
\sum_{i=1}^{n}( Y_j - \theta_0 ) \biggr|.
 \end{split}
\end{equation*}



But

$$\frac{1}{n} \,\biggr| \sum_{i=1}^{n}( Y_j - \theta_0 )
\biggr| = \bigr| \overline{Y} - \theta_0 \bigr| = \bigr|
\overline{Y} - E_Y \bigr| \rightarrow 0$$ in $G_{\theta_0, \eta_0}
-$probability, therefore Definition
\ref{CompositeDeconvolutionDef2} is satisfied if $\sqrt{n}
\,\bigr| \frac{1}{\sigma^2 (Y)} - \frac{1}{\widehat{\sigma}_n^2}
\bigr|$ is bounded in $G_{\theta_0, \eta_0} -$probability, and
this holds if $\widehat{\sigma}_n^2$ is a $\sqrt{n}-$consistent
estimate of $\sigma^2 (Y).$ Here $\overline{Y}$ denotes the sample
mean $\overline{Y} = \frac{1}{n} \,\sum_{i=1}^{n} Y_j .$
\end{proof}

\begin{proof} (The statement about $l^{\ast}_j$ from Example 3).

$$\frac{1}{\sqrt{n}} \,\biggr| \sum_{i=1}^{n}(l_j^\ast -
l^{\ast}_{\eta_0}(Y_j)) \biggr|
\quad\quad\quad\quad\quad\quad\quad\quad\quad\quad\quad\quad\quad
\quad\quad\quad\quad\quad\quad\quad\quad\quad\quad\quad$$

\begin{equation*}
 \begin{split}
 & = & \frac{1}{\sqrt{n}} \,\bigr| C_1(\eta_0) \bigr| \, \biggr|
\sum_{i=1}^{n}\bigr( (Y_j - \widehat{\theta}_n)^2 - (Y_j -
\theta_0)^2 \bigr) \biggr| \\
 & = & \frac{1}{\sqrt{n}} \,\bigr| C_1(\eta_0) \bigr| \, \biggr|
\sum_{i=1}^{n} (\widehat{\theta}_n - \theta_0) (-2 Y_j +
\widehat{\theta}_n + \theta_0) \biggr|
 \end{split}
\end{equation*}

\begin{equation*}
 \begin{split}
 & = & \bigr| C_1(\eta_0) \bigr| \,\, \sqrt{n\,} \bigr|
\widehat{\theta}_n - \theta_0 \bigr|
 \,\,\, \frac{1}{n} \,\biggr|
\sum_{i=1}^{n} ( Y_j - \widehat{\theta}_n ) + \sum_{i=1}^{n} ( Y_j
- \theta_0 ) \biggr| \\
 & = & \bigr| C_1(\eta_0) \bigr| \, \sqrt{n\,} \bigr|
\widehat{\theta}_n - \theta_0 \bigr|
 \, \,\bigr| ( \overline{Y} - \widehat{\theta}_n ) + (
\overline{Y} - \theta_0 ) \bigr| \\
 & \leq & \bigr| C_1(\eta_0) \bigr| \, \sqrt{n\,} \bigr|
\widehat{\theta}_n - \theta_0 \bigr|
 \,\,\Bigr(\bigr| \overline{Y} - \widehat{\theta}_n \bigr| + \bigr|
\overline{Y} - \theta_0 \bigr|\Bigr) \rightarrow 0
 \end{split}
\end{equation*}






\noindent in $G_{\theta_0, \eta_0} -$probability since for $n
\rightarrow \infty$ it holds that $\bigr| \overline{Y} -
\widehat{\theta}_n \bigr| \rightarrow 0$ and $\bigr| \overline{Y}
- \theta_0 \bigr| \rightarrow 0,$ both in $G_{\theta_0, \eta_0}
-$probability, and $\sqrt{n} \bigr| \widehat{\theta}_n - \theta_0
\bigr|$ is bounded in $G_{\theta_0, \eta_0} -$probability.
\end{proof}

\begin{proof}(Theorem \ref{CompositeDeconvolutionTh1}).
Put

\begin{equation}\label{CompositeDeconvolution16}
V_k =\biggr\{ \frac{1}{\sqrt{n}}\, \sum_{j=1}^{n}
l^{\ast}_{\theta_0}(Y_j)   \biggr\} \;
{(I^{\ast}_{\theta_0})}^{-1} \; {\biggr\{ \frac{1}{\sqrt{n}}\,
\sum_{j=1}^{n} l^{\ast}_{\theta_0}(Y_j) \biggr\} }^{T} \,,
\end{equation}

\noindent where $l^{\ast}_{\theta_0}$ is defined by
(\ref{CompositeDeconvolution11}) and $I^{\ast}_{\theta_0}$ by
(\ref{CompositeDeconvolution12}). Of course, $V_k$ is \emph{not} a
statistic since it depends on the unknown $\eta_0.$ But if the
true $\eta_0$ is known, then because of $\langle B1
\rangle$-$\langle B3 \rangle$ we can apply the multivariate
Central Limit Theorem and obtain $V_k \rightarrow_d
\,{\chi}^{2}_{k}$ as $n \rightarrow \infty\,.$ Condition
(\ref{CompositeDeconvolution14}) implies that
$$\frac{1}{\sqrt{n}} \sum_{i=1}^{n} l^\ast_j\,\, \rightarrow \,\,
\frac{1}{\sqrt{n}}\sum_{i=1}^{n}  l^{\ast}_{\theta_0}(Y_j)
\quad\mbox{in $G_{\theta_0, \eta_0}-$probability as}\,\, n
\rightarrow \infty$$ and by consistency of $\widehat{L}$ we get
the statement of the theorem by Slutsky's Lemma.
\end{proof}

\begin{proof}(Theorem \ref{CompositeDeconvolutionTh3}).
Because of condition (\ref{CompositeDeconvolution14}) the proof is
analogous to the proof of Theorem \ref{SimpleDeconvolutionTh8}.
Indeed, after obvious change of notations Propositions
\ref{SimpleDeconvolutionProp5}, \ref{SimpleDeconvolutionProp6},
and \ref{SimpleDeconvolutionProp7} are true for $W_k,$
$W_{S(l^\ast)},$ $S(l^\ast)$ instead of $U_k,$ $U_S,$ $S.$ Proofs
of the new versions of propositions are analogous to the proofs of
the previous versions. The only difference is that the proof of
the key inequality analogous to (\ref{SimpleDeconvolution23})
requires the use of the following lemma.

\begin{lemma}\label{CompositeDeconvolutionAlgebraicLemma1} Let $A$
be a $k \times k$ positive definite matrix and
$\{A_n\}^{\infty}_{n=1}$ be sequence of $k \times k$ matrices such
that $A_n\, \rightarrow \,A$ in the Euclidian matrix norm. Suppose
that for some real number $\delta > 0$ we have $A
>\delta$ in the sense that the matrix $(A- \delta I_{k \times k})$ is
positive definite, where $I_{k \times k}$ is the $k \times k$
identity matrix. Then for all sufficiently large $n$ it holds that
$A_n >\delta.$
\end{lemma}

\end{proof}
\end{document}